\documentclass[11pt]{article}
\usepackage{a4}
\usepackage{graphicx}
\usepackage{parskip}

\usepackage{amsfonts}
\usepackage{natbib}
\usepackage{algorithmic}

\usepackage{latexsym}

\newtheorem{theorem}{Theorem}[section]

\newtheorem{lemma}{Lemma}[section]
\newtheorem{corollary}{Corollary}[section]

\newtheorem{definition}{Definition}[section]

\newcommand{\R}{\mathbb{R}}
\newcommand{\Nat}{\mathbb{N}}

\newcommand{\PP} {{  \rm I\hskip-0.22em P}}

\newcommand{\EE} {{\rm I\hskip-0.48em E}}

\usepackage{a4}

\usepackage{graphicx}

\usepackage{parskip}

\usepackage{amsfonts}

\usepackage{natbib}

\usepackage{algorithm}

\usepackage{algorithmic}

\usepackage{amsmath, amssymb}

\usepackage{latexsym}

\begin{document}

\centerline{\bf Logistic regression with total variation regularization}

\vskip .1in

\centerline{Sara van de Geer, 5.3.2020}

\centerline{Seminar for Statistics, ETH Z\"urich} 

\vskip .1in

{\bf Abstract} We study logistic regression with total variation penalty on the canonical parameter and show
that the resulting estimator satisfies a sharp oracle inequality: the excess risk of the estimator is adaptive
to the number of jumps of the underlying signal or an approximation thereof. In particular
when there are finitely many jumps, and jumps up are sufficiently separated from jumps down, then the estimator
converges with a parametric rate up to a logarithmic term $\log n / n$, provided the tuning parameter
is chosen appropriately of order $1/ \sqrt n$. Our results extend earlier results for quadratic loss to logistic loss.
We do not assume any a priori known bounds
on the canonical parameter but instead only make use of the local curvature of the theoretical risk. 



\section{Introduction}\label{introduction.section}
In this paper we consider logistic regression with a total variation penalty on the canonical parameter.
Total variation based de-noising was introduced in \cite{rudin1992nonlinear}. Our aim here is to develop
theoretical results that show that the estimator adapts to the number of jumps in the signal.

For $i=1 , \ldots , n $, let $ Y_i\in \{0,1\}$ be independent binary observations. Write 
the unknown probability of success as $\theta_i^0 := P(Y_i=1) $,
and let $f_i^0 := \log ( \theta_i^0 / (1- \theta_i^0 ) )$ be the log-odds ratio, $i=1 , \ldots , n$. Define the total variation of
a vector $f \in \R^n$ as
$$ {\rm TV} (f) := \sum_{i=2}^n | f_i - f_{i-1} | . $$
We propose to estimate the unknown vector $f^0$ of log-odds ratios applying logistic regression with total variation
regularization. The estimator is
$$\hat f:= \arg \min_{f \in \R^n } \biggl \{ {1 \over n} \sum_{i=1}^n\biggl ( - Y_i f_i + \log (1+ {\rm e}^{f_i} ) \biggr ) + \lambda {\rm TV} (f) \biggr \} . $$

Our goal is to derive oracle inequalities for this estimator.
The approach we take shares some ideas with \cite{dalalyan2017prediction}, \cite{ortelli2018total} and \cite{vdGOrtelli2019}.
These papers deal with least squares loss, whereas the current paper 
studies logistic loss. Moreover, instead of using the projection arguments of the previous mentioned papers,
we use entropy bounds. This allows us to remove a redundant logarithmic term: we show that the 
excess risk of estimator $\hat f$ converges under certain conditions with rate $(s+1) \log n / n$ where $s$ is the number of jumps
of $f^0$ or of an oracle approximation thereof (see Theorem \ref{oracle.theorem}). This extends the result in
\cite{guntuboyina2020adaptive} - where there is also no redundant logarithmic term - to logistic loss and to a sharp
oracle inequality. 
 
To arrive at the results in this paper
we require that $\| \hat f\|_{\infty} $ stays bounded with high probability. In Theorem
\ref{bounded.theorem} we show that this requirement holds assuming that both $\| f^0 \|_{\infty}$ and
${\rm TV} (f^0) $ remain bounded. 

Theory for total variation regularization for least squares loss (the fused Lasso) has been developed in a series of papers (\cite{tibshirani2005sparsity}, \cite{tibshirani2014adaptive}, \cite{sadhanala2016total}, \cite{dalalyan2017prediction}, 
\cite{lin2017sharp}, \cite{padilla2017dfs}, \cite{Sadhanala2017}) including higher dimensional extensions
(\cite{hutter2016optimal}, 
\cite{chatterjee2019new}, \cite{fang2019multivariate}, \cite{Ortelli2d}) and higher order total variation
(\cite{steidl2006splines}, \cite{sadhanala2017higher}, \cite{vdGOrtelli2019}, \cite{guntuboyina2020adaptive}). 

Logistic regression with $\ell_1$-regularization has many applications. When there are co-variables, 
the penalty is on the total variation of the coefficients.
In \cite{yu2015classification} logistic regression with the fused Lasso is applied to spectral data, 
and \cite{liu2017structured} to gene expression data, whereas
\cite{ahmed2009recovering} applies it to time-varying networks.
In \cite{sun2012penalized} the penalty alternatively takes links between variables into account using a quadratic penalty.
The papers \cite{yu2015high} and \cite{liu2010efficient} present algorithms for fused Lasso.
In \cite{betancourt2017bayesian} a Bayesian approach with the fused Lasso is presented. 

This paper is organized as follows. In Section \ref{oracle.section} we state the oracle inequality for $\hat f$ 
(Theorem \ref{oracle.theorem}).
Section \ref{bounded.section} derives a bound for $\| \hat f \|_{\infty} $ (Theorem \ref{bounded.theorem}).
The remainder of the paper is devoted to proofs. Section \ref{standard.section} states some standard tools
to this end, Section \ref{oracle-proof.section} contains a proof of Theorem \ref{oracle.theorem} and 
Section \ref{Proof-bounded.section} a proof of Theorem \ref{bounded.theorem}. 


 \section{A sharp oracle inequality}\label{oracle.section}
 The empirical risk in this paper is given by the normalized minus log-likelihood
 $$ R_n (f) := {1 \over n} \sum_{i=1}^n\biggl ( - Y_i f_i + \log (1+ {\rm e}^{f_i} ) \biggr ), \ f \in \R^n .$$
 The theoretical risk is 
 $$R(f) := \EE R_n(f) , \ f \in \R^n  $$
 and $R(f) - R(f^0)$ is called the ``excess risk". For $f \in \R^n$, we write $\dot R_n (f) := \partial R(f) / \partial f$
 and $\dot R(f) := \EE \dot R_n (f) $. These are column vectors in $\R^n$. 
 Most of the arguments that follow go through for general convex differentiable loss functions.  
 We do use however that or all $f \in \R^n$, $\dot R_n (f) - \dot R(f) = -\epsilon^T f/n $ 
 where $\epsilon = Y - \EE Y $ is the noise. In other words, $f$ is the canonical parameter. In the
 case where the entries of the response vector $Y$ are in $\{0,1\}$, the entries of noise vector $\epsilon$ are 
 bounded by $1$. More generally, our theory would need that $\epsilon$ has sub-exponential entries.
 To avoid digressions, we simply restrict ourselves to logistic loss.
 
 Fix a vector ${\bf f} \in \R^n$. This vector will play the role of the ``oracle" as we will see in Theorem
 \ref{oracle.theorem}. We let 
$S:= \{ t_1 , \ldots , t_s \} $ ($1<t_1 < \cdots < t_s < n$) be the location
of its jumps: 
$${\bf f}_1 = \cdots = {\bf f}_{t_1-1} \not= {\bf f}_{t_1 } = \cdots = {\bf f}_{t_2-1} \not= {\bf f}_{t_2} \cdots
{\bf f}_{t_s-1} \not= {\bf f}_{t_s } = \cdots = {\bf f}_n . $$
Let $d_j := t_j - t_{j-1}$ be the distance between jumps, $j=1 , \ldots , r$,  where $r=s+1$, $t_r:=n+1$ and $t_0=1$.
Define $d_{\rm max} := \max_{1 \le j \le r } d_j $.

The quantities $\Delta_n^2$, $\delta_n^2 (t)$, $\lambda_n (t)$ and $\Gamma_n^2(t)$
we are about to introduce all depend on ${\bf f}$ although we do not express this
in our notation. Moreover, being non-asymptotic, these quantities are somewhat involved. After the
explicit expressions for $\Delta_n^2$, $\delta_n^2 (t)$ and $\lambda_n (t)$ we will give their
asymptotic order of magnitude. The asymptotic order of magnitude for $\Gamma_n^2 (t)$
depends on the situation. We discuss a special case after the statement of Theorem \ref{oracle.theorem}. 

We let 
 $$\Delta_n^2 := { 4  \sum_{j\in [1:r]:\ d_j \ge 1 }  (  \log (d_j-1) + 1) \over n} +  {s\over n}  ,$$
and define for $t>0$
\begin{eqnarray*}
 \delta_n^2 (t)  &:=&  \biggl ( {4 \nu  A_0 \Delta_n  }
 + 8 \sqrt { 1 + t +  \log (3+ 2\log_2 n)  \over n} \biggr )^2 \\
 &+ & \biggl ( {  2 \over  \nu }  
 +  4 \sqrt {A_0\Delta_n  \over n } 
 +      {4 \sqrt { 1 + t + \log (3+  2 \log_2 n)  }  \over n} \biggr )\\
 & & \times \biggl (  \Delta_n  + 2 \sqrt {s\over n}  \biggr )^2  , 
 \end{eqnarray*} 
 and 
 $$\lambda_n (t) := {1 \over \sqrt n} \biggl ( {  4 \over  \nu }  
 +  8 \sqrt {A_0 \Delta_n   \over  n } 
 +       {8 \sqrt { 1 + t + \log (3+  2\log_2 n)  }  \over n} \biggr ) . $$
 One sees that
  $$\Delta_n^2 = {\mathcal O} \biggl ({ (s+1) \log n \over n }  \biggr ) .$$
 Furthermore, for $\nu=1$ (say) and each fixed $t$ 
 $$\delta_n^2(t)= { \mathcal O}\biggl (  {(s+1)\log n \over n } \biggr ) , \ \lambda_n (t) = {\mathcal O}\biggl (   {1 \over \sqrt n} 
 \biggr ) ,$$
assuming $n^{-1} \sqrt {(s+1) \log n /n}   = {\mathcal O} (1)$ which is certainly
 true under the standard sparsity assumption $ {(s+1) \log n /n}= o(1)$. 

 The quantity $\delta_n^2 (t)$ will be part of the bound for the excess risk
 of $\hat f$, and $\lambda_n (t)$ can be thought of as the ``noise level" to be overruled by the penalty
 (see Theorem \ref{oracle.theorem}).  
 The constant $A_0$ is the (universal) constant appearing when bounding the entropy of the 
class of  functions  with both $\| \cdot\|_{\infty} $ and
${\rm TV} (\cdot )$ bounded by 1 (see Lemma \ref{monotone.lemma}).
The free parameter $t>0$ determines the confidence level of our statements. 
 Both $\delta_n (t)$ and $\lambda_n (t)$ depend on a further free parameter $\nu>0$
 which we do not express in our notation as one can simply choose
 $\nu=1$. It is however an option to choose $\nu$
 larger than 1, possibly growing with $n$: larger $\nu$ relaxes the requirement on the tuning parameter $\lambda$ but results in
 larger bounds for the excess risk. 
  
  Finally, we present a bound $\Gamma_n^2(t)$ for the so-called ``effective sparsity" as introduced in
 \cite{vdGOrtelli2019}, see also Definition \ref{Gamma.definition}. The effective sparsity may be seen
 as a substitute for the sparsity, which is defined as the number of active parameters of the oracle, which is $s+1$. The effective sparsity will in general be larger than $s+1$. Without going into details, we remark that this is due to correlations in the dictionary $X$ when writing
 $f = X b$, with dictionary $X \in \R^{n \times n}$ and coefficients $b_1:=f_1$, $b_k := f_k - f_{k-1}$, 
 $k \in [2:n ] $.  
 
 Let
 $q_{t_j} := {\rm sign} ({\bf f}_{t_j } )$, $j=1 , \ldots , s$. We write
  $J_{\rm monotone} :=  \{ 2\le j\le s:  \ q_{t_{j-1}} = q_{t_j} \} $ and 
  $J_{\rm change}:= [1:r] \backslash J_{\rm monotone}$. 
  Thus $J_{\rm monotone}$ are jumps with the same sign as the previous one, and
  $J_{\rm change}$ are jumps that change sign. We count the first jump as well as the endpoint $t_r=n+1$
  as  a sign change. 
Our bound for the effective sparsity is now
  $$ \Gamma_n^2 (t) := { \lambda_n^2 (t) \over  \lambda^2 }    \sum_{j \in J_{\rm monotone}  } { 8 (\log (d_j) +1 )  }+
  \sum_{j\in J_{\rm change} }{ 8n (\log (d_j ) + 2) \over d_j } 
  .$$
  
  The following theorem presents an oracle inequality for $\hat f$. Its proof can be found in
  Section \ref{oracle-proof.section}.
  
\begin{theorem} \label{oracle.theorem} 
Let ${\cal F}$ be a convex subset of $\R^n$  (possibly ${\cal F}= \R^n$) and 
$$ \hat f := \arg \min_{f \in {\cal F}} \biggl \{ R_n (f) + \lambda {\rm TV} (f) \biggr \} . $$
Assume  ${\bf f} \in {\cal F} $ satisfies $\| {\bf f}  \|_{\infty} \le B$  for some constant $B$ and define
$$\kappa  := { (1+ {\rm e}^B )^2 \over {\rm e}^B } . $$
Take 
$$\lambda \ge \lambda_n(t) \sqrt {d_{\rm max} \over 2 n} . $$
Then with probability at least $ \PP (\| \hat f \|_{\infty} \le  B) - \exp[-t] $ we have
$$ R(\hat f) - R({\bf f})\le 4 \kappa \delta_n^2 (t) + { \lambda^2 \over 4} \Gamma_n^2 (t). $$
\end{theorem}

Keeping the constant $B$ fixed, this theorem tells us that 
$$R( \hat f) - R({\bf f} ) = {\mathcal O}_{\PP} \biggl ( { \sum_{j=1}^r (\log (d_j ) +1 )  \over n }  + \lambda^2 \Gamma_n^2
\biggr ) $$
where we recall that $r=s+1$.
If the jumps of ${\bf f}$ are roughly equidistant, we see that $d_j \asymp d_{\rm max} \sim n/r$.
Taking $\lambda \asymp \lambda_n (t)/ \sqrt r \asymp \sqrt {1/(nr)}$, the bound for the effective sparsity
$\Gamma_n^2(t) $ is in the worst case (where the jumps of ${\bf f}$ have alternating signs)
of order $r^2 \log (n/r) $. In other words, in that case the rate is
$R( \hat f ) - R({\bf f} ) = {\mathcal O}_{\PP} (r \log (n/r) / n ) $, which for least squares loss is the minimax rate: see
\cite{lin2017sharp}.

If $\bf f$ is monotone, we get with $\lambda\asymp \sqrt {d_{\rm max} }/ n $
$$ \lambda^2 \Gamma_n^2 \asymp  {\sum_{j =2}^s   \log (d_j) +1 \over n } + 
{1 \over n}  \biggl ( { \log (d_1 ) d_{\rm max} \over d_1} + {\log (d_r ) d_{\rm max} \over d_r } \biggr ) .$$
In other words, the first jump of ${\bf f}$ should not occur to early, and the last jump not too late, relative
to the distance between the jumps. 

We note that the choice $\lambda \asymp \lambda_n (t) \sqrt{ d_{\rm max} / n }$ depends on
the oracle ${\bf f}$. Thus, if the tuning parameter $\lambda$ is given 
the choice of ${\bf f}$ 
depends on $\lambda$. 

We assumed that $\| {\bf f} \|_{\infty} \le B$. We do not assume $\| f^0 \|_{\infty} $ to be bounded by the
same constant $B$, 
but we do hope for a good approximation ${\bf f}$ of $f^0$ with $\| {\bf f} \|_{\infty} \le B$.
Nevertheless, Theorem \ref{oracle.theorem} presents a sharp oracle inequality
directly comparing $R(\hat f)$ with $R({\bf f})$: it does not require that the excess risk
$R({\bf f}) - R(f^0)$ is small in any sense. In the same spirit, the theorem
requires that $\| \hat f \|_{\infty}\le B$ with high probability. This can be accomplished by
taking ${\cal F} := \{ f \in \R^n  : \ \| f \|_{\infty} \le B \} $ (or some convex subset thereof).
Theorem \ref{oracle.theorem} holds for any $B$, i.e. it is a free parameter. 
However, one may not want to force
$\hat f$ to be bounded by a given constant but let the data decide for a bound on $\hat f$. This is a reason
why we establish Theorem \ref{bounded.theorem} given in the next section.

\section{Showing that $\| \hat f \|_{\infty} $ is bounded (instead of assuming this)}\label{bounded.section}
Since $f^0$ minimizes $R(f)$
a two-term Taylor expansion around $f^0$ gives
$$ R(f) - R(f^0) ={1 \over 2} \ddot R (\tilde f)$$
where $\tilde f_i$ lies between $f_i $ and $f_i^0$, $i=1 , \ldots , n $.
It follows that
$$R(f) - R(f_0) \ge {1 \over 2 K_f^2 }  \| f - f^0 \|_{Q_n}^2  $$
where 
$$\| \cdot \|_{Q_n} = \| \cdot \|_2/ \sqrt n $$
and where (for logistic loss)
$${ K_f^2 } :=  { (1+
{\rm e}^{\| f \|_{\infty} \vee \| f^0 \|_{\infty} } )^2 \over {\rm e}^{\| f \|_{\infty} \vee \|f^0\|_{\infty}}}  .$$
Thus, if both $\| f \|_{\infty} $ and $ \|  f^0\|_{\infty}$ stay within bounds we have standard quadratic curvature
of $R (\cdot)$ at $f^0$. 
Otherwise the
the constant $K_f$ grows
exponentially fast. We will therefore assume that $\| f^0 \|_{\infty}$ stays bounded and our task is
then to show that $\| \hat f  \|_{\infty}$ stays bounded as well.
The following theorem (where we have not been very careful with the constants) is derived in Section
\ref{Proof-bounded.section}.

\begin{theorem}\label{bounded.theorem} Let ${\rm TV} (f^0) \le M_0 $ for some constant
$M_0 \ge 1$. Define
$$ K:=  {(1 + {\rm e}^{1 + 2^4 M_0 + \| f^0 \|_{\infty}  } )^2 \over {\rm e}^{1 +2^4  M_0 + \| f^0 \|_{\infty} } } . $$
Suppose 
\begin{eqnarray*} 
\lambda &\le &  \biggl ( 2^4 (2K^2) M_0 \biggr )^{-1} \\
 \lambda &\ge&  2^{8}  n^{-2/3} A_0^{2/3} (2K^2)^{1/3}  \\
 \lambda  &\ge & 2^{8} (2K^2)  { 1+t \over n} .
\end{eqnarray*} 
where the last inequality holds for some $t >0$, and where in the second last inequality $A_0$ is the constant appearing when bounding the entropy of the 
class of  functions  with both $\| \cdot\|_{\infty} $ and
${\rm TV} (\cdot )$ bounded by 1 (see Lemma \ref{monotone.lemma}).
  Then with probability at least $1- \exp[-t] $ it holds that
  $$ { \| \hat f - f^0 \|_{Q_n}^2 \over 2K^2 } + \lambda {\rm TV} (\hat f- f^0) \le
  4 \lambda M_0  $$
 and
  $$ \| \hat f - f^0 \|_{\infty}  \le {1 +8M_0 \over 2}. $$
  \end{theorem}
  
  One may object that the conditions on the tuning parameter $\lambda$ depend on $f^0$ via bounds
  on $\| f^0 \|_{\infty}$ and ${\rm TV} (f^0)$. On the other hand, the choice of $\lambda$ in
  Theorem \ref{oracle.theorem} will be of larger order than $n^{-2/3}$ if one
  aims at adaptive results, and it will need to tend to zero. For such $\lambda$ and
  for $\| f^0 \|_{\infty}$ and ${\rm TV} (f^0)$ remaining bounded, the conditions of
  Theorem \ref{bounded.theorem} will be met for all $n$ sufficiently large.

\section{Some standard results useful for both Theorem \ref{oracle.theorem} and
Theorem \ref{bounded.theorem}} \label{standard.section}

 \begin{lemma}\label{Hoeffding.lemma}
 We have for all vectors $g \in \R^n$, 
 $$\PP ( \epsilon^T g \ge \| g \|_2 \sqrt {2t} ) \le \exp [-t] , \ \forall \ t > 0 . $$
  \end{lemma}
  {\bf Proof.}
  The entries in $\epsilon$ have mean zero, are bounded by 1,
and are independent. This means we can apply Hoeffding's inequality to $\epsilon^T g$.
\hfill $\sqcup \mkern -12mu \sqcap$

For ${\bf Q}$ a  probability measure on $\{ 1, \ldots , n \}$ and 
a set ${\cal G}\subset \R^n$  we let $H(\cdot , {\cal G} , {\bf Q})$ be
the entropy\footnote{For $u>0$ the $u$-covering number $N(u)$ of a metric space $({\cal V} , d)$ is the smallest $N$
such that there exists $\{ v_j \}_{j=1}^N \subset {\cal V} $ with
$\sup_{v \in {\cal V} } \min_{1 \le j\le N} d (v , v_j) \le u $. The entropy is $H( \cdot) := \log N( \cdot )$.}
of ${\cal G}$ endowed with the metric induced by the $L_2 ({\bf Q})$-norm

\begin{lemma} \label{entropy-integral.lemma} Let ${\cal G} \subset \R^n $ be a set with
diameter
$$R:= \sup_{g \in {\cal G} } \| g \|_{Q_n}  . $$
Suppose
$$J( R):= 2 \int_0^R \sqrt {2 H( u , {\cal G} , Q_n)} du $$ 
exists. Then for all $t >0$, with probability at least $1- \exp[-t]$
it holds that
$$\sup_{g \in {\cal G} } \epsilon^T g /n \le { J (R) \over \sqrt n} +  4 R \sqrt {1+t \over  n} . $$
\end{lemma} 

{\bf Proof.} 
We can apply Hoeffding's inequality to $\epsilon^T g$ for
each $g$ fixed, see Lemma \ref{Hoeffding.lemma}.
The result of the current lemma is thus essentially applying Dudley's entropy integral. The constants are taken from Theorem
17.3 in \cite{vdG2016}.
\hfill $\sqcup \mkern -12mu \sqcap$

\begin{lemma} \label{monotone.lemma}
Let ${\cal G} := \{ g \in \R^n : \  \| g \|_{\infty} \le 1 , \ {\rm TV} ( g ) \le 1 \} $. 
It holds for any probability measure ${\bf Q}$ 
$$ H( u, {\cal G} , {\bf Q} ) \le{ A_0 \over u } \ \forall \  u > 0 $$
where $A_0$ is a universal constant. 
\end{lemma} 

{\bf Proof.} See \cite{van1996weak}, Theorem 2.7.5. \hfill $\sqcup \mkern -12mu \sqcap$

 \section{Proof of Theorem \ref{oracle.theorem}.} \label{oracle-proof.section}

\subsection{The main body of the proof of Theorem \ref{oracle.theorem}.} \label{main0body.section}
The following lemma is Lemma 7.1 in \cite{vdG2016}. We present a proof for completeness.
\begin{lemma} \label{KKT.lemma} Let ${\cal F}$ be a convex subset of $\R^n$ (possibly ${\cal F} = \R^n$) and
$$ \hat f := \arg \min_{f \in {\cal F}} \biggl \{ R_n (f) + \lambda {\rm TV} (f) \biggr \} . $$
Then for all $f \in {\cal F}$
$$- \dot R_n (\hat f)^T (f- \hat f) \le \lambda{\rm TV} (f) - \lambda {\rm TV} (\hat f) . $$
\end{lemma}
{\bf Proof of Lemma \ref{KKT.lemma}.}
Define for $0< \alpha < 1$, $\hat f_{\alpha} := (1- \alpha) \hat f + \alpha f $. Then,  using the convexity of
${\cal F} $ 
\begin{eqnarray*}
R_n (\hat f) + \lambda{\rm TV} (\hat f)& \le & R_n (\hat f_{\alpha} ) + \lambda{\rm TV} (\hat f_{\alpha} ) \\
& =&  R_n (\hat f_{\alpha} ) + (1- \alpha) \lambda{\rm TV} (\hat f) + \alpha \lambda {\rm TV} (f) .
\end{eqnarray*} 
Thus
$$ { R_n (\hat f) - R_n (\hat f_{\alpha} ) \over \alpha } \le \lambda{\rm TV} (f) - \lambda{\rm TV} (\hat f) . $$
The result now follows by letting $\alpha \downarrow 0$.
\hfill $\sqcup \mkern -12mu \sqcap$

\begin{lemma} \label{basic-inequality.lemma} 
Let ${\cal F}$ be a convex subset of $\R^n$ and
$$ \hat f := \arg \min_{f \in {\cal F}} \biggl \{ R_n (f) + \lambda {\rm TV} (f) \biggr \} . $$
Then for all $f \in {\cal F}$ 
$$ R(\hat f) - R(f) + {\rm rem} (f , \hat f )  \le \epsilon^T ( \hat f - f)/n + \lambda {\rm TV} (f) - \lambda {\rm TV} (\hat f ) , $$
where
$$ {\rm rem} (f, \hat f ) = R(f) - R(\hat f) - \dot R (\hat f )^T ( f - \hat f ) . $$
\end{lemma}

{\bf Proof of Lemma \ref{basic-inequality.lemma}.}
By Lemma \ref{KKT.lemma} 
$$- \dot R_n (\hat f)^T (f- \hat f) \le \lambda{\rm TV} (f) - \lambda {\rm TV} (\hat f) . $$
So
\begin{eqnarray*}
 R(\hat f) - R(f) + {\rm rem} (f , \hat f ) & =&  - \dot R (\hat f)^T ( f - \hat f ) \\
& =& (\dot R_n (\hat f) -\dot R (\hat f ))^T (f - \hat f ) - \dot R_n (\hat f )^T (f - \hat f) \\
& = & \epsilon^T ( \hat f - f) / n - \dot R_n (\hat f )^T (f - \hat f) \\
& \le & \epsilon^T ( \hat f - f)/n + \lambda {\rm TV} (f) - \lambda {\rm TV} (\hat f ) . 
\end{eqnarray*} 
\hfill $\sqcup \mkern -12mu \sqcap$

One sees from Lemma \ref{basic-inequality.lemma} that we need appropriate bounds
for the empirical process $\{ \epsilon^T f: \ f \in \R^n\}$. These will be established in the next
two subsections, Subsections \ref{EPT.section} and \ref{materialEPT.section}.
In Subsection \ref{EPT.section} we announce the final result, and Subsection
\ref{materialEPT.section} presents the technicalities that lead to this result. 

\subsection{The empirical process $\{ \epsilon^T f: \ f \in \R^n \} $} \label{EPT.section} 

We 
 consider the weights\footnote{These weights are inspired by the following. Let ${\cal V}_S$ be the linear space of functions that are piecewise constant with jumps at $S$
and $\Pi_S$ be the projection operator on the space ${\cal V}_S$. Then
$$\epsilon^T f / n = \epsilon^T \Pi_{S}f / n + \epsilon^T (I- \Pi_S) f / n ,$$
and one can verify that
$$ \epsilon^T (I- \Pi_Sf ) / n = \sum_{j \notin S} V_k (f_k - f_{k-1} ) $$
where $V_{-S} = \{ V_k \}_{k \notin S} $ is a vector of random variables
with ${\rm var} ( V_k)= w_k^2$, $k \notin S$. }
$$w_{k}^2  := \begin{cases} \biggl ( { k- t_{j-1}  \over d_j } \biggr  ) \biggl ({ t_j - k \over n } \biggr ) , & t_{j-1} +1 \le 
k \le t_{j} -1,  \ j \in [1:r]\cr \ \ \ \ \ \ \ {1 \over n},& k= t_{j}, \ j\in [1:s] \cr  \end{cases}   . $$
For a vector $f\in \R^n$ we define $(Df)_k:= f_k - f_{k-1}$ ($k=[2:n]$) so that
$\| Df \|_1= {\rm TV} (f)$. 
Let $w= (w_1 , \ldots , w_n )^T $ be the vector of weights and $w^{-1} := ( 1/w_1 , \ldots , 1/ w_n)$.
Write
$$ w_{-S} (Df)_{-S} := \{ w_k (Df)_k \}_{k \notin S } . $$
We use the notation $\| \cdot \|_{Q_n} := \| \cdot \|_2/ \sqrt n$ for the normalized
Euclidean norm. 
For $t >0$ let 
\begin{eqnarray*}\label{delta.equation}
 \delta_n^2 (t)&\ge&  \biggl ( {4 \nu  A_0\| w^{-1} \|_{Q_n} \over  \sqrt n }
 + 8 \sqrt { 1 + t +  \log (3+ 2\log_2 n)  \over n} \biggr )^2 \\
 &+ & \biggl ( {  1 \over 2 \nu }  
 +  4 \sqrt {A_0\| w^{-1} \|_{Q_n}/ \sqrt n  \over n } 
 +      {4 \sqrt { 1 + t + \log (3+  2 \log_2 n)  }  \over n} \biggr )\\
 & & \times \biggl (  \| D w \|_2 + 2 \sqrt {s\over n}  \biggr )^2  , 
 \end{eqnarray*} 
 and 
 $$\lambda_n (t)  \ge   {1 \over \sqrt n} \biggl ( {  4 \over  \nu }  
 +  8 \sqrt {A_0 \| w^{-1} \|_{Q_n} / \sqrt n  \over  n } 
 +       {8 \sqrt { 1 + t + \log (3+  2\log_2 n)  }  \over n} \biggr ) . $$

After establishing  the material of Subsection \ref{materialEPT.section} we are able show the following result:
\begin{theorem}\label{EPT.theorem} Let $\mu>0$ and $t>0$ be arbitrary.
With probability at least $1- \exp[-t]$
 $$ \epsilon^T f / n \le \mu \delta_n^2(t) + {\| f \|_{Q_n} ^2 \over \mu } + \lambda_n (t) \| w_{-S} D_{-S}  f  \|_1  $$
 uniformly for all $f \in \R^n$. 
\end{theorem}

{\bf Proof of Theorem \ref{EPT.theorem}.} This follows from combining
Lemma \ref{partial-integration.lemma}  with Lemma \ref{EPT.lemma}
(see Corollary \ref{EPT.corollary}). 
\hfill $\sqcup \mkern -12mu \sqcap$ 

\subsection{Material for the result for the empirical process $\{ \epsilon^T f: \ f \in \R^n \} $ in Theorem \ref{EPT.theorem}}
\label{materialEPT.section}

Let for all $f \in \R^n$, 
$$\gamma_f: =
{ \sum_{k=1}^n f_j / w_j \over \| w^{-1} \|_2^2 } $$
and let 
$$f_{\rm P} :=\Pi_{w^{-1}}  f := w^{-1} \gamma_f $$
be the projection of $f $ on the vector $w^{-1} $. 
Define the anti-projection $f_{\rm A} := (I- \pi_{w^{-1}} ) f$. 

We let 
$$wf:= \{ w_k f_k \}_{k=1}^n .$$

We start with some preliminary bounds.

\begin{lemma} \label{preliminary-bounds.lemma} For all $f \in \R^n$ it holds that 
$$ \| wf - \gamma_f \|_{\infty} \le {\rm TV} (wf) $$ 
and
$$ { \| f_{\rm A} \|_{\infty}\over {\rm TV} (wf) } \le \sqrt n . $$
\end{lemma} 
{\bf Proof of Lemma \ref{preliminary-bounds.lemma}.} 
For all $i \in [1:n] $,
\begin{eqnarray*}
w_i f_i - \gamma_f  &=& w_i f_i - { \sum_{k=1} f_k / w_k \over \| w^{-1} \|_2^2 } \\
& = & { \sum_{k=1}^n (w_i f_i - w_k f_k  ) / w_k^2 \over \| w^{-1} \|_2^2 } \le {\rm TV} (wf) . 
\end{eqnarray*}
or $\| wf - \gamma_f \|_{\infty} \le {\rm TV} (wf)$.
Since when $g = wf$
$$f_{\rm A} = w^{-1} (g- \gamma_f ) $$
we see that
$$ \| f_{\rm A} \|_{\infty} \le \|w^{-1} \|_{\infty} {\rm TV} (g) = \| w_{-1} \|_{\infty} {\rm TV } (wf ) . $$
Since $\| w_{-1} \|_{\infty} = \sqrt n $ we conclude that
$$ \| f_{\rm A} \|_{\infty} \le \sqrt n {\rm TV} (wf) . $$
\hfill $\sqcup \mkern -12mu \sqcap$

We use Dudley's entropy integral to bound the empirical process
over $\{ f: \ \| f_{\rm A} \|_{Q_n} \le R, \ {\rm TV} (wf) \le 1\}$ with the radius $R$ some fixed value.

\begin{lemma} \label{fixedR.lemma} Let $R >0$ be arbitrary. For all $t >0$, with probability at least $1- \exp[-t]$, 
\begin {eqnarray*}
\sup_{ \| f_{\rm A}  \|_{Q_n}  \le R , \ {\rm TV} (wf) \le 1 } \epsilon^T f/n  &\le &
   4 \sqrt {2A_0\| w^{-1} \|_{Q_n}   R \over n}   + 4 R \sqrt {1+t \over n}   .
\end{eqnarray*}
\end{lemma}

{\bf Proof of Lemma \ref{fixedR.lemma}.}  
Let ${\bf Q}_w$ be the discrete probability measure
that puts  mass 
$ {w_i^{-2} / \| w^{-1}  \|_2^2}    $ on $i$, ($i \in [1: n ] $). 
Denote
the $L_2({\bf Q}_w)$-norm by $\| \cdot \|_{{\bf Q}_w} $. 
For ${\cal G} \subset \R^n$ we let
${\cal H} ( \cdot , {\cal G}, {\bf Q}_w   )$ denote the entropy of ${\cal G}$ for the metric induced by $\| \cdot \|_{{\bf Q}_w} $.
By Lemma \ref{preliminary-bounds.lemma}
$$  \| wf - \gamma_f \|_{\infty}  \le {\rm TV} (wf)  . $$
Thus by Lemma \ref{monotone.lemma}, with $A_0$ the constant given there,
$${\cal H} ( u , \{ wf - \gamma_f  :\  \ {\rm TV} (wf)  \le 1 \} ,{\bf Q}_w   ) \le {A_0 \over u } \ \forall \ u>0 .$$
For $f \in \R^n$ we have
$$ \| f_{\rm A}  \|_{Q_n}^2 = {1 \over n} \sum_{i=1}^n ( w_i f_i - \gamma_f )^2 / w_i^2 =
\|  wf - \gamma_f \|_{{\bf Q}_w}^2    \| w^{-1}  \|_{Q_n}^2     .$$
Therefore
$${\cal H} ( u , \{ f_{\rm A} ,  \  {\rm TV} (wf ) \le 1  \} , Q_n  ) 
\le {A_0\| w^{-1} \|_{Q_n}   \over u}   \ \forall \ u>0 . 
$$
The entropy integral can therefore be bounded as follows
$$ 2 \int_0^R  \sqrt {2 {\cal H}  ( u , \{ f_{\rm A}  :  \| f_{\rm A}  \|_{Q_n}   \le  R , \  {\rm TV} (wf ) \le 1  \} , Q_n ) } du$$
$$ \le4 \sqrt {2 A_0   \| w^{-1} \|_{Q_n}  R }  . $$
By Lemma \ref{entropy-integral.lemma}
the result follows.
\hfill $\sqcup \mkern -12mu \sqcap$

The next lemma invokes Lemma \ref{fixedR.lemma} and the peeling device to
obtain a result for the weighted empirical process. 

\begin{lemma} \label{varyingR.lemma}
For all $t>0$, with probability at least $1- \exp[-t]$ it holds that
\begin{eqnarray*}
\epsilon^T f_{\rm A} /n &\le&  8 \sqrt {A_0 \| w^{-1} \|_{Q_n}   \over n} \biggl (  \sqrt { \| f_{\rm A} \|_{Q_n} 
{\rm TV} (wf) }\vee { {\rm TV} (wf) \over  n^{3/4}}  \biggr )  \\
& + & 8 \biggl ( \| f_{\rm A}  \|_{Q_n}   \vee { {\rm TV} (wf) \over  n^{3/2} }  \biggr )  \sqrt { 1 + t + \log (2+  2\log_2 n )   \over n}  
\end{eqnarray*}
uniformly over all $f$.
\end{lemma}

{\bf Proof of Lemma \ref{varyingR.lemma}.} Let $t>0$ and let ${\cal A} $ be the event
\begin{eqnarray*}
\biggl \{ \epsilon^T f_{\rm A} /n &\ge &  8 \sqrt {A_0\| w^{-1} \|_{Q_n} \over n } \sqrt { \| f \|_{Q_n} \vee {1 \over n^{3/2}} }
\\
& + & 8 \biggl (\| f  \|_{Q_n}  \vee {1 \over n^{3/2}} \biggr )  \sqrt { 1 + t +  \log(2+  2 \log_2 n  ) \over n}  , \\
\mbox{for some} \ &f&  \mbox {with} \   { \| f   \|_{Q_n}  } \le \sqrt n\ \mbox{and}\ 
{\rm TV} (wf) \le 1 \biggr \} .
\end{eqnarray*}
Let ${\cal A}_0$ be the event
\begin{eqnarray*}
 \biggl \{ \sup_{   \| f_{\rm A} \|_{Q_n}  \le {1 \over n^{3/2}}   , \ {\rm TV} (wf ) \le 1 }
\epsilon^T f_{\rm A} /n &\le& 8 \sqrt {A_0\| w^{-1} \|_{Q_n} \over  n}  \sqrt {1 \over n^{3/2}} \\
& +& {8 \over n^{3/2}}   \sqrt {1+t  + \log (2+2 \log_2n )   \over n}  \biggr \}  .
\end{eqnarray*}
Let $N \in \Nat$ satisfy $2 \log_2 n \le N \le 1+ 2 \log_2 n $ and for $j \in [1: N] $ let ${\cal A}_j$ be the event
\begin{eqnarray*}
 \biggl \{ \sup_{  {2^{j-1} \over n^{3/2} } < \| f_{\rm A} \|_{Q_n} \le { 2^{j} \over n^{3/2}}  , \ {\rm TV} (wf ) \le 1 }
\epsilon^T f_{\rm A} /n &\le& 8 \sqrt {A_0\| w^{-1} \|_{Q_n } \over  n}  \sqrt {2^{j-1} \over n^{3/2}} \\
& +& { 8 2^{j-1} \over  n^{3/2} } \sqrt {1+t  + \log (2+2\log_2 n )  \over n}  \biggr \}  .
\end{eqnarray*}
Application
of  Lemma \ref{fixedR.lemma} gives that for all $j \ge 0$, 
$$ \PP ({\cal A}_j ) \le \exp[-(t+ \log (2 + 2 \log_2 n  ) ] . $$
Since ${\cal A}  \subset \cup_{j=0}^N {\cal A}_j $ follows that
$$ \PP ({\cal A} ) \le \sum_{j=0}^N \PP ({\cal A}_j) \le (1+N)  \exp[-(t+\log (2+2 \log_2 n )] \le \exp[-t] . $$
The result now follows by replacing
$f_{\rm A}$ by
$f_{\rm A} / {\rm TV } (wf) $ and noting that
$$ {\rm TV} \biggl (w f_{\rm A} / {\rm} TV (wf)  \biggr ) =1,$$
and invoking from Lemma \ref{preliminary-bounds.lemma} the bound
$$ \| f_{\rm A} / {\rm TV} (wf) \|_{Q_n} \le \| f_{\rm A} / {\rm TV} (wf) \|_{\infty}  \le \sqrt n .$$
\hfill $\sqcup \mkern -12mu \sqcap$

We present a corollary that applies the ``conjugate inequality"
$2ab \le a^2 + b^2 $ (with constants $a$ and $b$ in $\R$), then gathers terms and applies the conjugate inequality again. 

\begin{corollary} \label{varyingR.corollary} Let $\nu>0$ and $\mu>0$  be arbitrary.
For all $t>0$ with probability at least $1- \exp[-t]$
\begin{eqnarray*}
& & \epsilon^T f_{\rm A} /n\\ &\le&  \biggl ( { 4\nu  A_0\| w^{-1} \|_{Q_n}  \over \sqrt n} 
 + 8 \sqrt { 1 + t +  \log (2+ 2\log_2 n)  \over n} \biggr ) {  \| f_{\rm A} \|_{Q_n}}   \\
 &+&\biggl ( {  4 \over  \nu }  
 +  8 \sqrt {A_0\| w^{-1} \|_{Q_n}/ \sqrt n  \over n } 
 +      {8 \sqrt { 1 + t + \log (2+  2\log_2 n)  }  \over n} \biggr ){ {\rm TV} (wf) \over \sqrt n}  \\
&\le& {\mu \over 2}  \biggl ( {4\nu  A_0\| w^{-1} \|_{Q_n}  \over \sqrt n} 
 + 8 \sqrt { 1 + t +  \log (2+ 2 \log_2 n)  \over n} \biggr )^2 \\ & +& { \| f_{\rm A} \|_{Q_n}^2 \over 2 \mu}  \\
 &+&\biggl ( {  4 \over  \nu }  
 +  8 \sqrt {A_0 \| w^{-1} \|_{Q_n}/ \sqrt n  / \over n  } 
 +      {8 \sqrt { 1 + t + \log (2+  2 \log_2 n)  }  \over n} \biggr ) { {\rm TV} (wf) \over  \sqrt n } 
\end{eqnarray*}
uniformly for all $f$. 
\end{corollary}

We now add the missing $f_{\rm P} = f - f_{\rm A}$.

\begin{lemma}\label{EPT.lemma}
For all $t>0$ with probability at least $1- \exp[-t]$
\begin{eqnarray*}
& & \epsilon^T f /n\\ 
&\le& {\mu \over 2 }  \biggl ( {4  \nu  A_0\| w^{-1} \|_{Q_n} \over \sqrt n} 
 + 8 \sqrt { 1 + t +  \log (3+ 2 \log_2 n)  \over n} \biggr )^2 \\ & +& { \| f \|_{Q_n}^2 \over 2 \mu}  \\
 &+&\biggl ( {  4 \over  \nu }  
 +  8 \sqrt {A_0 \| w^{-1} \|_{Q_n} / \sqrt n  \over n  } 
 +      {8 \sqrt { 1 + t + \log (3+  2 \log_2 n)  }  \over n} \biggr ) { {\rm TV} (wf) \over  \sqrt n } 
\end{eqnarray*}
uniformly for all $f$. 
\end{lemma}

{\bf Proof of Lemma \ref{EPT.lemma}.}
We have by Pythagoras' rule $ \| f \|_2^2 = \| f_{\rm P} \|_2^2 + \| f_{\rm A} \|_2^2$.
Moreover, by Hoeffding's inequality, with probability at least $1- \exp[-t]$
$$ \epsilon^T f_{\rm P} / n \le \| f_{\rm P} \|_{Q_n}  \sqrt {2 t \over n} \le
{ \mu t \over n} + { \| f_{\rm P} \|_{Q_n}^2  \over 2 \mu }. $$
\hfill $\sqcup \mkern -12mu \sqcap$

In Lemma \ref{EPT.lemma} the term including ${\rm TV} (wf)$ is almost but not yet quite the one to be dealt with
by the penalty. We bound it by $\| w_{-S} D_{-S} f \|_1$ with appropriate remaining terms
invoking the ``chain rule".  Here
$$ w_{-S} (Df)_{-S} := \{ w_k (Dk)_k \}_{k \notin S } .$$

\begin{lemma}\label{partial-integration.lemma} For all $f \in \R^n$
$${\rm TV} (wf) \le  \sqrt n \biggl ( \| Dw \|_2+ 2 \sqrt {s / n } \biggr )\| f \|_{Q_n}  +
\| w_{-S} D_{-S} f \|_1.$$
\end{lemma}

{\bf Proof of Lemma \ref{partial-integration.lemma}.} 
We use that 
\begin{eqnarray*}
 {\rm TV} (wf)  &\le& \sum_{i=2}^n |(w_i - w_{i-1} ) f_{i-1} | +\sum_{i=2}^n | w_{i} ( f_i - f_{i-1} ) |  \\
& \le & \| D w \|_2   \| f \|_2  + \| w D f  \|_1 .
\end{eqnarray*} 
Moreover
$$\| w D f \|_1 = \| w_S D_S f \|_1 + \| w_{-S} D_{-S} f \|_1 $$
with
$$ w_{S} (Df)_{S} := \{ w_k (Dk)_k \}_{k \in S } ,$$
satisfying
\begin{eqnarray*}
 \| w_S D_S f \|_1 &=& \sum_{j=1}^{s} |f_{t_{j}+1} - f_{t_j}| / \sqrt n \\
 & \le &  \sqrt s  
\sqrt {\sum_{j=1}^s |f_{t_{j}+1} - f_{t_j}|^2} /\sqrt n \\
&\le & 2 \sqrt s  \| f \|_2 /\sqrt n .
\end{eqnarray*}
Thus
\begin{eqnarray*} 
{\rm TV} (wf)& \le&  \biggl ( \| Dw \|_2+ 2 \sqrt {s / n } \biggr )\| f \|_2  +
\| w_{-S} D_{-S} f \|_1.
\end{eqnarray*} 
\hfill $\sqcup \mkern -12mu \sqcap$

\begin{corollary} \label{EPT.corollary} The result from Theorem \ref{EPT.theorem} 
now follows using
$$ \biggl (\| D w \|_2 + 2 \sqrt {s /n} \biggr ) \| f \|_{Q_n} \le
{\mu \over 2} \biggl (\| D w \|_2 + 2 \sqrt {s /n} \biggr )^2 + {\| f \|_{Q_n}^2 \over 2 \mu } , \ f \in \R^n . $$

\end{corollary} 

\subsection{Bounds for the weights and their inverses}\label{weights.section}

So far we assumed in this section (see Subsection \ref{EPT.section}), 
that for $t>0$, the quantities $\delta_n^2 (t)$ and $\lambda_n (t)$ involved in the bound for the empirical process in Theorem \ref{EPT.theorem}
satisfy
\begin{eqnarray*}
& & \delta_n^2 (t) \ge   \biggl ( {4 \nu  A_0\| w^{-1} \|_{Q_n} \over  \sqrt n }
 + 8 \sqrt { 1 + t +  \log (3+ 2\log_2 n)  \over n} \biggr )^2 \\
 &+ & \biggl ( {  1 \over 2 \nu }  
 +  4 \sqrt {A_0\| w^{-1} \|_{Q_n}/ \sqrt n  \over n } 
 +      {4 \sqrt { 1 + t + \log (3+  2 \log_2 n)  }  \over n} \biggr )\\
 & & \times \biggl (  \| D w \|_2 + 2 \sqrt {s\over n}  \biggr )^2  , 
 \end{eqnarray*} 
 and 
 $$\lambda_n (t) \ge  {1 \over \sqrt n} \biggl ( {  4 \over  \nu }  
 +  8 \sqrt {A_0 \| w^{-1} \|_{Q_n} / \sqrt n  \over  n } 
 +       {8 \sqrt { 1 + t + \log (3+  2\log_2 n)  }  \over n} \biggr ) . $$
involving $\| w^{-1} \|_{Q_n} $ and $\| Dw \|_2 $.
In this subsection, we present bounds for these, so leading to the values
$\delta_n^2 (t)$ and $\lambda_n(t) $ presented in Section \ref{oracle.section}. 

\begin{lemma}\label{weights.lemma} It holds that
 $$ \| w^{-1} \|_2^2 \le 2 n \sum_{d_j \ge 2} (  \log (d_j-1) + 1) + n s \le n^2 \Delta_n^2 $$
 and
 $$ \| Dw \|_2^2 \le 4 \sum_{d_j \ge 2}  (\log (d_j-1) + 1) /n+ s/n=:  \Delta_n^2 .$$
\end{lemma}

{\bf Proof of Lemma \ref{weights.lemma}.} 
We have\footnote{We use $\sum_{k=1}^{d-1} { d \over k (d- k)} =
\sum_{k=1}^{d-1} \biggl ( {1 \over k} + { 1 \over d-k} \biggr ) = 2 \sum_{k=1}^{d-1} {1 \over k} \le 
2 (1+\log (d-1) )$.} 
\begin{eqnarray*}
 \| w^{-1} \|_2^2 &=& \sum_{d_j \ge 2}  \sum_{k=1}^{d_j-1 } { n d_j \over k ( d_j-k)} + {ns } \\
& \le & 2 n \sum_{j=1}^r(\log (d_j-1) + 1) + ns 
\end{eqnarray*}
Moreover, for $1\le k \le d_j-1 $, $j \in [1:r] $, 
$$ |\sqrt k \sqrt {d_j-k} - {\sqrt {k-1} \sqrt {d_j-(k-1)} } |
\le \sqrt {d_j - k  \over k} + \sqrt {k-1 \over d_j - k } $$
$$ \le \sqrt {d_j -1\over  k} + \sqrt  {d_j-2 \over d_j- k} \le \sqrt {d_j \over  k} + \sqrt  {d_j \over d_j- k}  $$
so that 
\begin{eqnarray*}
& &  \sum_{k=1}^{d_j -1} { |\sqrt k \sqrt {d_j-k} - {\sqrt {k-1} \sqrt {d_j-(k-1)} } |^2 \over n d_j } \\ &\le&
{ 2 \over n} 
\sum_{k=1}^{d_j-1} \biggl ( {1 \over  k} + { 1 \over d_j-k  } \biggr ) \\
& \le& {1 \over n} \sum_{j=1}^r (4\log (d_j-1) + 2) . 
\end{eqnarray*} 
Finally, for $j \in [1:s]$
$$ | w_{t_j } - w_{t_j -1} | = \biggl | {1 \over \sqrt n} - \sqrt { d_j-1 \over d_j} { 1 \over \sqrt n} \biggr |
 \le {1 \over \sqrt n} .$$
\hfill $\sqcup \mkern -12mu \sqcap$

\subsection{A bound for the effective sparsity}\label{effective-sparsity.section}

We let for all $f \in \R^n$
$$(Df)_S := \{ (Df)_k \}_{k \in S} , \ (Df)_{-S} := \{ (Df)_k \}_{k \notin S } . $$
and recall that
$$w_{-S} (Df)_{-S} := \{ w_k (Dk)_k \}_{k \notin S } . $$
Let $q_{t_j}:= {\rm sign} ({\bf f}_{t_j}) $, $j\in [1:s] $. We define
$q_S := \{ q_{t_j} \}_{j=1}^s $.

\begin{definition} \label{Gamma.definition} Let $\lambda \ge \lambda_n (t) \sqrt { d_{\rm max} / (2n) } $.
The effective sparsity at ${\bf f}$ is
$$ \Gamma^2 ({\bf f},t):=\left (\min \biggl \{ \| f \|_{Q_n} : \ q^T (Df )_S - \| ( 1- w_{-S} \lambda (t) / \lambda ) (D f)_{-S}  \|_1 =1 \biggr \}\right )^{-2} . $$
\end{definition} 

Recall the definitions
$$J_{\rm monotone} := \{ 2\le j \le s : q_{t_j} = q_{t_{j-1} } )\}  , \ J_{\rm change} := [1:r]/ J_{\rm monotone} .$$

\begin{lemma}\label{Gamma.lemma}
For 
$\lambda \ge \lambda_n (t) \sqrt {d_{\rm max} / n}  $ we have
$$ \Gamma({\bf f} , t) \le \Gamma_n^2 (t) , $$
where
\begin{eqnarray*}
\Gamma_n^2 (t ) &:= & { \lambda_n^2 (t) \over  \lambda^2 }  \sum_{j \in J_{\rm monotone}   } { 8 (\log (d_j) +1 )  }
+ \sum_{j\in J_{\rm change} } { 8n (\log (d_j ) + 2) \over d_j } 
   . 
\end{eqnarray*} 
\end{lemma}

{\bf Proof of Lemma \ref{Gamma.lemma}.}
The proof uses interpolating vectors $q \in \R^{n} $ as in \cite{vdGOrtelli2019} where $q= (q_1 , q_{-1} )^T$ is given below.
We show that
$$ q_S^T (Df)_S - \|  ( 1- w_{-S} \lambda (t) / \lambda ) (D f)_{-S}  \|_1 \le
q_{-1}^T D ( {\bf f } - \hat f ) .$$
The result then follows from
$$ q_{-1}^T D ( {\bf f } - f ) = (D^T q_{-1})^T ( {\bf f } - f) \le \| D^T q_{-1} \|_2 \| {\bf f } - f \|_2 .$$
Furthermore, under the boundary conditions $q_1 = q_n= 0$ we see
that $\| D^T q_{-1} \|_2 = \| D q \|_2 $. 
 Define
$$\omega_k^2 := \begin{cases} \biggl ( { k- t_{j-1}  \over d_j } \biggr  ) \biggl ({ t_j - k \over n } \biggr ) 
{\lambda_n (t) \over \lambda } & t_{j-1} +1 \le 
k \le t_{j} -1,  \ j \in J_{\rm monotone}, \ d_j \ge 2 \cr
 \biggl ( { k- t_{j-1}  \over d_j } \biggr  ) \biggl ({ t_j - k \over d_j } \biggr ) & t_{j-1} +1 \le 
k \le t_{j} -1,  \ j \in J_{\rm change} \cr \ \ \ \ \ \ \ 0 & k= t_{j}, \ j\in [1:s] \cr  \end{cases}   $$
For $j \in [1:r]$ we let $\bar t_j = {t_{j-1} + t_j \over 2}$ be the midpoints.
Moreover, for $k \notin \{ t_1 , \ldots , t_s\}$ let 
$$q_k := \begin{cases} 0 & 1 \le k < {\bar t_1}\cr  {\rm sign}( {\bf f}_{t_1} ) (1-2 \omega_k ) &   \bar t_1  \le 
k \le t_{1} -1 \cr 
{\rm sign}( {\bf f}_{t_{j-1}} ) (1-2 \omega_k ) &    t_{j-1} +1 \le 
k < \bar t_{j}, \ j \in [2 :s]\cr 
{\rm sign} ({\bf f}_{t_{j}} ) (1-2 \omega_k ) &  \  \bar t_j  \le 
k \le t_{j} -1, \ j \in [2 :s]  \cr 
{\rm sign} ({\bf f}_{t_{r-1}} ) (1-2 \omega_k ) &  \  t_{r-1}    \le 
k < \bar t_r  \cr
0 & \bar t_r \le k \le n \cr 
\end{cases} $$

We get that for $ \bar t_j -1 \le k < \bar t_j$, $j \in J_1$
$$ |1- 2\omega_k| \le { 4 \over d_j } $$

For $j \in J_{\rm monotone} $ we see that
$$ \sum_{k=1}^{d_j}   |q_{t_{j-1}+k } - q_{t_{j-1} + k-1} |^2 \le {\lambda_n^2 (t) \over \lambda^2 } { 8 ( \log d_j +1 ) \over  n} 
  . $$
and for $j \in {J_{\rm change}}$, 
$$ \sum_{k=1}^{d_j} |q_{t_{j-1}+k } - q_{t_{j-1} + k-1} |^2  \le { 8 (\log d_j + 2)  \over d_j} . $$
Thus
$$ \| D q \|_2^2 \le { \lambda_n^2 (t) \over  \lambda^2 }  \sum_{j \in J_{\rm monotone}   } { 8 (\log (d_j) +1 ) \over n  } 
+ \sum_{j\in J_{\rm change} } { 8 (\log (d_j ) + 2) \over d_j } 
  . $$
  The lemma now follows from $\Gamma^2 ({\bf f} , t) \le n \| D q \|_2^2 $. 
\hfill $\sqcup \mkern -12mu \sqcap$

\subsection {Finalizing the proof of Theorem \ref{oracle.theorem}} \label{finalizing,section}
We have by Lemma \ref{basic-inequality.lemma}
\begin{eqnarray*}
 & & R(\hat f) - R({\bf f}) + {\rm rem} ({\bf f} , \hat f )\\
   & \le& 
\mu \delta_n^2  (t)+ {\| \hat f - {\bf f}  \|_{Q_n}  \over \mu} + \lambda_n (t) \| w_{-S}  D_{-S}\hat f \|_1 +
\lambda \| D_S {\bf f} \|_1 - \lambda \|D \hat f \|_1 \\
& =&  \mu \delta_n^2  (t)+ {\| \hat f - {\bf f} \|_{Q_n}^2 \over \mu} +  \lambda \biggl (   \| D_S {\bf f} \|_1 -\| D_S \hat f \|_1
- \| (1- \lambda_n (t) w_{-S} / \lambda )D_{-S} \hat f \|_1\biggr )  \\
& \le & \mu \delta_n^2  (t)+ {\| \hat f - {\bf f}  \|_{Q_n}^2 \over \mu} + \lambda \Gamma_n (t)   \| \hat f -
{\bf f} \|_{Q_n} \\
& \le & \mu \delta_n^2 (t)  + {2 \| \hat f - {\bf f}  \|_{Q_n}^2  \over \mu} + {\lambda^2 \over 4} \Gamma_n^2 (t). 
\end{eqnarray*} 
Choose $\mu= 4 \kappa$ to obtain
$${2 \| \hat f - {\bf f}  \|_{Q_n}^2  \over \mu} = { \| \hat f - {\bf f}  \|_{Q_n}^2 \over 2 \kappa} \le
{\rm rem} ({\bf f} , \hat f ) $$
whenever $\|\hat f \|_{\infty} \le B$.
\hfill $\sqcup \mkern -12mu \sqcap$

\section{Proof of Theorem \ref{bounded.theorem}} \label{Proof-bounded.section}

\subsection{Some lemmas used in the proof of Theorem \ref{bounded.theorem}} \label{lemmas.section}
The proof of Theorem \ref{bounded.theorem} applies some auxiliary lemmas which we develop
in this subsection.
Define
$$ \tau (f ) := \| f \|_{Q_n}/(\sqrt 2K) + (\lambda /\delta ) {\rm TV} (f) $$
with
$$\delta^2 := 2^4 \lambda M_0 , \ 
K^2:= { (1+
{\rm e}^{ 1+2^4 M_0 + \| f^0 \|_{\infty} } )^2 \over {\rm e}^{1+ 2^4 M_0 + \|f^0\|_{\infty}}}  $$
where
$M_0 \ge {\rm TV} (f^0) \vee 1$. Moreover, we let
$$ \hat f_{\alpha} := \alpha \hat f + (1- \alpha ) f^0 .$$
with $$ \alpha := {\delta  \over \delta  + \tau (f-f_0)}  .$$
Let ${\cal F}_0 := \{ f :\ \tau(f ) \le \delta \} $.

\begin{lemma} \label{delta.lemma}
It holds that $\hat f_{\alpha} - f^0 \in {\cal F}_0$, i.e., 
$ \tau( \hat f_{\alpha} - f^0 )  \le \delta $.
Moreover, if in fact 
$ \tau( \hat f_{\alpha} - f^0 )  \le \delta/2 $, 
then also $\hat f- f^0 \in {\cal F}_0 $. 
\end{lemma}

{\bf Proof.} 
We have
$$ \tau ( \hat f_{\alpha} - f^0 ) = \alpha \tau(\hat f - f^0) =
{ \delta \tau(\hat f - f^0) \over \delta  + \tau (\hat f-f_0)} \le \delta . $$
If in fact $ \tau( \hat f_{\alpha} - f^0 )  \le \delta/2 $, 
we have
$$\tau ( \hat f_{\alpha} - f^0 ) = { \delta \tau(\hat f - f^0) \over \delta  + \tau (\hat f-f_0)} \le \delta/2 $$
which gives 
$ \tau(\hat f - f^0) \le \delta /2  + \tau (\hat f-f_0)/2$
or 
$ \tau( \hat f - f^0) \le \delta $.
\hfill $\sqcup \mkern -12mu \sqcap$

\begin{lemma} \label{average.lemma} For all $f \in \R^n$
$$ \| f \|_{\infty} \le \| f \|_{Q_n}  + {\rm TV} (f) . $$
Moreover, 
$${\cal F}_0 \subset \{ f : \ \| f \|_{\infty} \le \sqrt 2K \delta +  \delta^2 / \lambda, \ {\rm TV} (f) \le \delta^2 /\lambda \} . $$
\end{lemma}

{\bf Proof.} For $f \in \R^n$ we denote its average by
$$ \bar f := {1 \over n} \sum_{i=1}^n f_i .$$
Then
$$ \| f \|_{Q_n}^2 = \bar f^2 +\| f - \bar f \|_{Q_n} \ge \bar f^2   .$$
Moreover, for all $i$,
$$ f_i - \bar f = {1 \over n} \sum_{j=1}^n ( f_i - f_j ) \le {\rm TV} (f) . $$
It follows that
$$ \| f \|_{\infty} \le \bar f + \| f - \bar f \|_{\infty} \le \| f \|_{Q_n}   + {\rm TV} (f) . $$
For $f \in {\cal F}_0$ we have $\| f \|_2/\sqrt n \le \sqrt 2K\delta $ and
${\rm TV}  (f) \le \delta^2/ \lambda $ so that also
$\| f \|_{\infty} \le \sqrt 2 K \delta + \delta^2 /\lambda $.
\hfill $\sqcup \mkern -12mu \sqcap$

\begin{lemma} \label{K.lemma} 
Let
$$K^2 :={ (1+
{\rm e}^{ 1+2^4M_0 + \| f^0 \|_{\infty} } )^2 \over {\rm e}^{1+ 2^4 M_0 + \|f^0\|_{\infty}}} $$
and let $\delta^2 := 2^4  \lambda M_0 \le  1/( 2K^2)$.
Then for all $f$ with $f-f^0 \in {\cal F}_0$ it is true that $K_f \le K$.
\end{lemma}

{\bf Proof.} Since 
for $f - f^0 \in {\cal F}_0$, $\| f - f^0 \|_{\infty} \le \sqrt 2 K \delta + \delta^2 / \lambda \le 1+ 2^4 M_0$, we see that
$\| f \|_{\infty} \le 1+ 2^4 M_0 + \| f ^0\|_{\infty} $.
Therefore
$$ K_f^2 = { (1+
{\rm e}^{\| f \|_{\infty} \vee \| f^0 \|_{\infty}  } )^2 \over {\rm e}^{\| f \|_{\infty} \vee\| f^0 \|_{\infty} } } \le K^2.$$
\hfill $\sqcup \mkern -12mu \sqcap$

 \begin{lemma} \label{basic-inequality.lemma} We have
 $$ 0 \le R(\hat f ) - R( f^0)\le \epsilon^T ( \hat f - f^0)/ n
 + \lambda {\rm TV} (f^0) - \lambda {\rm TV} (\hat f) .$$
 This inequality is also true with $\hat f$ replaced by $\hat f_{\alpha}$.
 \end{lemma}
 
 {\bf Proof.} 
 For any $f$
 \begin{eqnarray*}
 0 \le R(f ) - R( f^0) & =&
 - \biggl [ \biggl ( R_n (f ) - R(f ) \biggr )- \biggl ( R_n (f^0) - R(f^0) \biggr ) \biggr ] \\
 &+& R_n ( f ) - R_n (f^0) \\
 &=& \epsilon^T (f- f^0) / n + R_n (f) - R_n (f^0) .
 \end{eqnarray*} 
  Insert  the basic inequality
 $$ R_n (\hat f) + \lambda {\rm TV} (\hat f) \le R_n (f^0) + \lambda {\rm TV} (f^0)  $$
 or 
 $$R_n(\hat f ) - R_n (f^0) \le \lambda {\rm TV} (f^0 ) - \lambda {\rm TV} (\hat f )  $$
 to arrive at the first statement of the lemma.
 To obtain the second statement, we note that
 by convexity of $f \mapsto R_n (f)$ such basic inequality is also true for $\hat f_{\alpha}$:
 \begin{eqnarray*}
& &  R_n (\hat f_{\alpha} )+ \lambda {\rm TV} (\hat f_{\alpha} ) \\
& \le & \alpha R_n (\hat f) + \alpha \lambda {\rm TV} (\hat f) + (1- \alpha)  R_n (f^0) + (1- \alpha) \lambda {\rm TV} (f^0)\\
& \le & R_n (f^0) + \lambda {\rm TV} (f^0) . 
 \end{eqnarray*}
 \hfill $\sqcup \mkern -12mu \sqcap$

%

\subsection{Proof of Theorem \ref{bounded.theorem}} \label{bounded.section}

We have for $f \in {\cal F}_0 $, $\| f \|_{\infty} \le \sqrt 2K \delta + \delta^2 / \lambda \le 2 \delta^2 / \lambda$
and as well as ${\rm TV} (f) \le \delta^2 / \lambda \le 2 \delta^2 / \lambda $.
It follows from Lemma \ref{monotone.lemma} that
$$H (u,  {\cal F}_0  , Q_n ) \le {2A_0\delta^2 \over \lambda u } \ \forall \ u >0 $$
so that
\begin{eqnarray*}
2 \int_0^{\sqrt 2 K \delta} \sqrt { 2 H(u, {\cal F}_0  , Q_n ) } du
&\le& 4 \sqrt {2A_0\sqrt 2 K } {\delta \over \sqrt \lambda} \int_0^{\sqrt 2 K \delta } { 1 \over \sqrt u } du  \\
&=& 8 \sqrt {2 A_0  \sqrt 2 K \over \lambda } \delta^{3/2} . 
\end{eqnarray*}
But then, in view of Lemma \ref{entropy-integral.lemma}, for all $t>0$ with probability at least $1-\exp [-t] $,
\begin{eqnarray*}
 \sup_{f \in {\cal F}_0  } \epsilon^T f/n &\le& 
8 \sqrt {2 A_0\sqrt 2 K \over n \lambda } \delta^{3/2}  + 4\sqrt 2 K  \delta \sqrt { 1+ t\over n } .
\end{eqnarray*}
Since, by Lemma \ref{delta.lemma}, $ \hat f_{\alpha } - f^0  \in {\cal F}_0 $
we know from Lemma \ref{K.lemma} that
$K_{\hat f_{\alpha} } \le  K $.
Thus, in view of Lemma \ref{basic-inequality.lemma} and the
bound 
$$R( \hat f_{\alpha} )- R(f^0) \ge { \| \hat f_{\alpha} - f^0 \|_{Q_n}^2  \over 2 K^2} ,$$
 we have shown that with probability at least $1- \exp[-t]$
\begin{eqnarray*}
 & & { \| \hat f_{\alpha} - f^0 \|_{Q_n}^2 \over 2K^2 } + 
\lambda {\rm TV} ( \hat f_{\alpha }- f^0 )\\ & \le&  2 \lambda {\rm TV } (f^0) +
 8 \sqrt {2 A_0\sqrt 2 K \over n \lambda } \delta^{3/2} + 4 \sqrt 2 K \delta \sqrt { 1+ t\over n } \\
& \le&  \  2\lambda M_0 \ \  +\ \ 
 8 \sqrt {2 A_0\sqrt 2 K \over n \lambda } \delta^{3/2}  + 4 \sqrt 2 K \delta \sqrt { 1+ t\over n }.
\end{eqnarray*}
We want the three terms on the right hand side to add up to at most $\delta^2 / 4$.
We choose
\begin{eqnarray*}
 \lambda M_0 &=& \delta^2 / 2^3 \\
 8 \sqrt {2 A_0\sqrt 2 K \over n \lambda } \delta^{3/2}  &\le& \delta^2 / 2^4, \\
4 \sqrt 2 K \delta \sqrt { 1+ t \over n } &\le& \delta^2 / 2^4 . 
\end{eqnarray*}
or
\begin{eqnarray*}
 { 2^4 \lambda M_0 } &=& \delta^2  \\
 \biggl ( { 2^7  \sqrt {2 A_0 \sqrt 2 K}  \over \sqrt {n  \lambda}   } \biggr )^4&\le& \delta^2 , \\
\biggl ( {2^6  \sqrt 2 K  }  \sqrt { 1+ t \over n   }  \biggr )^2 &\le & \delta^2  . 
\end{eqnarray*}
The first one is the largest of the three.
This leads to the requirements 
$$ { 2^4 \lambda M_0 } \ge 
 \biggl (  2^7 \sqrt {2 A_0\sqrt 2 K \over n \lambda }  \biggr )^4   $$
 which is true for 
  $$ \lambda \ge 2^{8}  n^{-2/3} A_0^{2/3} (\sqrt 2K)^{2/3}  $$
  and
$$ 2^4  \lambda M_0 \ge  \biggl ( {2^6  }  \sqrt 2 K \sqrt { 1+ t \over n   }  \biggr )^2   $$
which holds for
$$\lambda  \ge 2^{8} (2K^2)  { 1+t \over n} $$
  where we invoked for both requirements that $M_0 \ge 1$. 
  Then with probability at least $1- \exp[-t]$
  $$ { \| \hat f_{\alpha} - f^0 \|_{Q_n}^2  \over 2 K^2}  +\lambda {\rm TV} (\hat f_{\alpha}- f^0 ) \le \delta^2 / 4 . $$
    For all $f \in \R^n$
    $$ \delta \tau ( f) =  {\delta \| f \|_{Q_n} \over \sqrt 2 K } + \lambda {\rm TV} (f) \le
  \delta^2 / 4 +  {  \| f \|_2^2 / n \over 2 K^2 } + \lambda {\rm TV} (f) . $$
   Thus we have shown that
  $$\delta \tau (\hat f_{\alpha} - f^0) \le \delta^2/ 4 + \delta^2 / 4 = \delta^2 / 2 $$
  or
  $$ \tau (\hat f_{\alpha} - f^0) \le \delta / 2 . $$
  By Lemma \ref{delta.lemma} this implies $\hat f \in {\cal F}_0$. We can now apply the same
  arguments to $\hat f$ as we did for $\hat f_{\alpha}$ to obtain that with probability at
  least $1- 2 \exp[-t] $ it holds that
  $$ { \| \hat f - f^0 \|_{Q_n}^2 \over 2K^2 } + \lambda {\rm TV} (\hat f- f^0) \le
  \delta^2 / 4 = 4 \lambda M_0 . $$
  This implies by Lemma \ref{average.lemma}
  $$ \| \hat f - f^0 \|_{\infty} \le {\sqrt 2 K \delta\over 2} + {\delta^2 \over 4 \lambda } \le {1 +8M_0 \over 2}. $$
   
  \hfill $\sqcup \mkern -12mu \sqcap$

\bibliographystyle{plainnat}
\bibliography{reference}
\end{document}